\documentclass{amsart}

\usepackage{amssymb, amscd, amsmath, amsthm, 
epsf, epsfig, 
latexsym}

\newtheorem{theorem}{Theorem}

\newtheorem{observation}[theorem]{Observation}

\theoremstyle{definition}

\theoremstyle{remark}

\numberwithin{equation}{section}

\begin{document}

\title{Observations on Lickorish Knotting of Contractible 4--Manifolds}

\author{Charles Livingston}

\address{Department of Mathematics, Indiana University, Bloomington, IN 47405}

\email{livingst@indiana.edu}

\keywords{\texttt{Contractible 4--Manifolds}}

\subjclass{57Q45}

\date{\today}

\begin{abstract} Lickorish has constructed large families of contractible
4--manifolds that have knotted embeddings in the 4--sphere and has also shown that
every finitely presented perfect group with balanced presentation occurs as the
fundamental group of the complement of  a knotted contractible manifold. 
Here we make a few observations regarding Lickorish's construction, showing how to
extend it to construct contractible 4--manifolds which have an
infinite number of knotted embeddings and also to construct knotted embeddings of
the Mazur manifold for which the complement has trivial fundamental group.
\end{abstract}

\maketitle

\vskip.2in 
In his recent paper, Lickorish \cite{li} describes a clever construction yielding
for each finitely presented perfect group $G$ with balanced presentation a compact
contractible 4--manifold $M_G$ with two embeddings in $S^4$, one for which the
complement is diffeomorphic to $M_G$ and the other with complement having
fundamental group $G$.  Here we make several observations based on Lickorish's
work.  With minor modifications, we use the notation of \cite{li} throughout.  

Thanks go to Ray Lickorish, first for identifying this interesting topic,
secondly for pointing out some of the issues resolved here, and finally for his
reflections on my  first attempts at extending his results.

\vskip.2in

\begin{observation} The construction can be modified to assure
that:
(1) For each group $G$ there is an infinite family of $M_G$ having the
desired  pair of embeddings and (2) For different groups $G$ the
constructed infinite  families have no elements in common.
\end{observation}

\begin{proof} We let $M_0$ be  the  Mazur manifold \cite{ma} with Kirby
diagram as shown in Figure 1.  (The curves $\gamma$ and $\gamma'$ are extraneous
for now.)  Since 
$M_0$ embeds in
$S^4$ with simply connected complement (as in \cite{li}, $ M_0 \times I \cong
B^5$), the  manifold
$M_G$ of the construction can be replaced with the boundary connected sum
$M_G \#_{\partial} M_0$.  This manifold will still have two embedding into $S^4$,
one with contractible complement and the other with complement having
fundamental group $G$.  Note that this manifold is not diffeomorphic to
$M$
since the boundary has changed by forming the connected sum with the
boundary of
the Mazur manifold, which Mazur  showed is not  $S^3$.  By repeating this
process
one can easily build   the desired families of examples; for instance,
arrange that for the first group $G$ all the boundaries have a prime
number of summands, for the second group $G$ arrange that all have a
prime squared number of summands, etcetera.
\end{proof}

\begin{figure}[t]
\begin{center}
\epsfxsize=2in
\epsfbox[94 158 310 348]{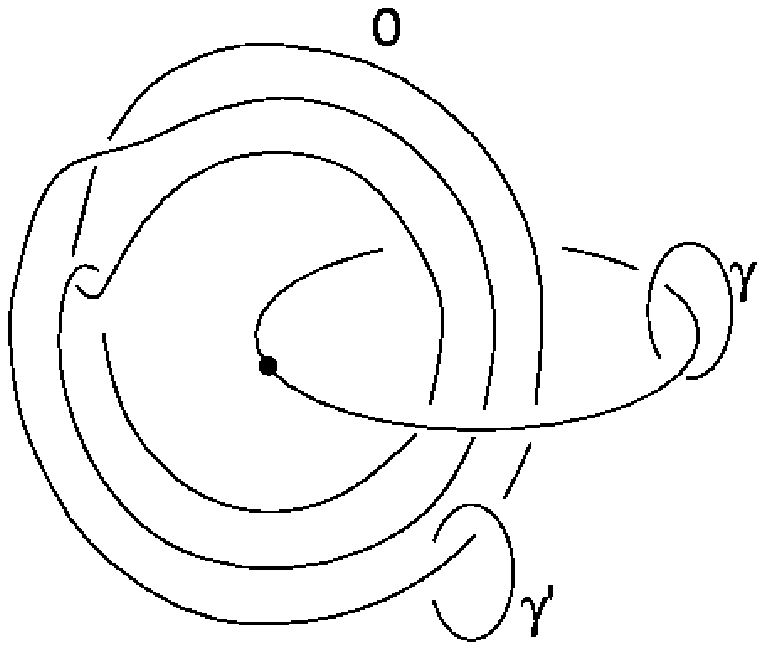}
\caption{}
\end{center}
\end{figure}

Following Observation \ref{obs2} we will give examples showing that groups with 
the
desired properties exist. 
\vskip.1in

\begin{observation}\label{obs2}
If the finitely generated  perfect group $G$  with balanced presentation maps
onto a nontrivial finite quotient of a 2--knot group then
there
exists an infinite family of embeddings, $\{\phi_i\}$, of the manifold constructed
by Lickorish,
$M_G$,  into
$S^4$ distinguished by
$\pi_1(S^4 - \phi_i(M_G))$.
\end{observation}
\vskip.1in
\begin{proof} Let $K$ be the 2--knot and let $H = \pi_1(S^4 -
K)/N$
be
the finite quotient.   Let $\rho \! : G \to H$ be the given surjective
homomorphism. Also, let $x$ be a meridian of $K$ and let $\bar{x}$ denote
its image in $H$.

Pick an element $g \in G$ such that $\rho(g) = \bar{x}$.  It can be
arranged in
the construction of $M_G$ that $g$ is among the generators of $G$ in the
balanced
presentation: simply add a generator $z$ to the balanced presentation
along with
the
 relation $z = g$, with $g$ written in terms of the generators of the
initial
balanced presentation.

It now follows from the initial construction of $M_G$ as the complement of
$X_G$
that
the meridian of the 2--handle of $M_G$ corresponding to the generator $z$
represents $g \in G =\pi_1(X_G)$.
To construct a new embedding of $M_G$ into $S^4$, tie the knot $K$ as a
local knot
in the given
2--handle of $M_G$.  This does not change $M_G$ but changes the fundamental
group of
$X_G$; the new group is constructed from the free product $G * \pi_1(S^4 -
K)$ by
identifying $g \in G$ with $x \in \pi_1(S^4 - K)$. We denote this group by
$G_1$ and also write it as $G *_Z \pi_1(S^4 -
K)$ though it is not an amalgamated product in the case that $g$ has 
finite order in
$G$. (It is not clear that this new group is not isomorphic to $G$.)

There are two homomorphisms of $\pi_1(S^4 - K)$ to $H$: the first is the
projection, $p$, the second factors through the cyclic group $Z$, mapping
the
meridian to $\bar{x}$.  Denote this second map by $q$.  The maps $\rho *
p$ and  $\rho * q$ each determine homomorphisms of $G_1$ to $H$.  These
homomorphisms
are
distinct as one is surjective when restricted to the image of $\pi_1(S^4 -
K)$
and the other is not surjective when restricted to this subgroup.  (Note
that
$H$
is perfect since it is a quotient of $G$ and hence is not cyclic.)

We first observe that these two embeddings cannot be isotopic; if they  were
 there would be an isomorphism from  $  G$ to  $G *_Z \pi_1(S^4 -
K)$ carrying the meridian   representing to $g$ to the meridian $m'$ representing
$g = x$.   Thus there would be a group
isomorphism from $G$ to $G *_Z
\pi_1(S^4 - K)$ sending $g$ to $g = x$.  However that cannot be as, by the above
argument, $G$ and $G *_Z
\pi_1(S^4 - K)$ have different numbers of homomorphisms onto $H$ sending these 
preferred meridians to $\bar{x}$.  (Notice that since $H$ is finite there is a
finite number of such homomorphisms.)

 By
repeating the construction of locally knotting the 2--handle of $M$ using
$K$, a 
sequence of nonisotopic embeddings of $M$ into $S^4$ is constructed.

We cannot show that the sequence of fundamental groups of the complements are all
distinct, but by counting homomorphisms to $H$ it follows that some subsequence 
must
be distinct; the number of homomorphisms onto $H$ goes to infinity since after
adding $n$ knots to the band there are at least 
$2^n$ homomorphisms onto
$H$. 

\end{proof}
\vskip.2in
\noindent{\bf Examples\ }  This simplest example of   Observation 2
occurs
with the binary icosahedral group, $H(2,3,5)$, the perfect group with 120
elements
representing the fundamental group of +1 surgery on the trefoil knot.
This
group
clearly has a balanced presentation and is also a quotient of the trefoil
group,
which is isomorphic to the fundamental group of the 0--twist spin of the
trefoil.  
More generally, consider the group of the $r$--fold cyclic branched cover
of the
$(p,q)$--torus knot, denoted $H(p,q,r)$.  If $p$, $q$, and $r$, are
pairwise
relatively prime, then $H(p,q,r)$ is perfect.  Furthermore, according to
Gordon,
\cite{go}, the $r$--twist spin of the $(p,q)$--torus knot has fundamental group
$H(p,q,r) \times Z$, and hence maps onto $H(p,q,r)$.  It remains to show
that
$H(p,q,r)$ has a finite quotient.  A presentation of $H(p,q,r)$ is given
by
$<x,y| x^p = y^q = (xy)^r>$.  For   such groups a nontrivial
representation to an alternating group can be constructed.  (This was
done by Fox in \cite{fo}; a more accessible  reference is Milnor \cite{mi}.)

\vskip.2in 
\begin{observation}
  There exist contractible manifolds   built with a single 1--handle that possess
two embeddings in $S^4$, one with simply connected complement and one with
nontrivial  complementary fundamental group. 
\end{observation}

\vskip.1in
\begin{proof} This  fact follows from the result of Neuzil \cite{ne} showing
that the Dunce Cap embeds in $S^4$ with nonsimply connected complement.  The
following approach gives us more control over the contractible manifold as well. 
Suppose that 
$L = L_1
\cup L_2$ is a 2--component link in $S^3$ with $L_1$ unknotted bounding a trivial
slice disk $D_1$ in $B^4$ and $L_2$   slice, bounding a slice disk $D_2$ in $B^4$.
Assume the linking number is 1. Let $S^3$ separate $S^4$ into two components, 
$B_1$
and
$B_2$ and view $D_1 \subset B_1$ and $D_2 \subset B_2$.  

As in Lickorish's construction, we let $M = (B_1 - N(D_1)) \cup N(D_2)$.  This is
clearly a contractible 4--manifold that doubles to give $S^4$.  However, its
complement is $X   = (B_2 - N(D_2)) \cup N(D_1)$.  Its group is given by the group
of the complement of the slice disk with an added relations coming from adding the
2--handle, $N(D_1)$.
\end{proof}
\vskip.1in

\noindent{\bf Examples\ }  For any knot $K$, the knot $L_2 = K \# -K$ is slice 
with
fundamental group of the complement of the slice disk being $\pi_1(S^3 - K)$.  Any
element of this group can be represented by an unknot $L_1$ in the complement of
$L_2$.  Hence, the groups that arise in this construction include all perfect
groups constructed by adding a single relation to a classical knot group.  For 
instance:  {\sl the fundamental group of a
homology 3--sphere built by surgery on a classical knot is the fundamental group of
the complement of the embedding of a contractible 4--manifold with one 1--handle 
into
$S^4$.}
 
\vskip.2in

In the previous construction it is not clear that we are constructing distinct
embeddings if $L_2$ is unknotted; this case is perhaps the most perplexing.  We
have the following example.
\vskip.1in
\begin{observation}
The Mazur Manifold illustrated in Figure 1 has two nonisotopic embeddings into
$S^4$.
\end{observation}
\vskip.1in
\begin{proof}  As described for instance in  \cite{ak}, the boundary of the
Mazur manifold
$M$ has an involution $F$ carrying $\gamma$ to $\gamma'$.  A handlebody picture of
the  manifold
$M
\cup_{id}  M$ is formed from Figure 1 by adding a 2--handle with 0 framing to
$\gamma'$.  Similarly, a handlebody picture of the 
manifold
$M
\cup_{F}  M$ is formed from Figure 1 by adding a 2--handle with 0 framing to
$\gamma$.  (In each case a 3-- and 4--handle must be added as well.)  It is an easy
exercise in Kirby calculus \cite{aki} to see that both are $S^4$.  Hence, we have
two embeddings of $M$ into $S^4$.

Clearly, in the first case the curve $\gamma'$ is slice in the complement --- the
2--handle is added to $\gamma'$.  In the second case $\gamma'$ is not slice in the
complement --- this is a result of Akbulut, \cite{ak}. 
\end{proof}

\vskip.2in

\noindent{\bf Questions}  In the above construction, if $L_2$ is unknotted  is
there an embedding of the constructed Mazur--like manifold into $S^4$ with
nonsimply connected complement?  If the crossing that is not part of the clasp  of
the attaching map of the 2--handle in Figure 1 is changed the previous argument 
does
not apply --- Akbulut has shown that in this case $\gamma$ will be slice.  Does this
manifold knot in $S^4$?  Does there exist a contractible 4--manifold that does not
knot in $S^4$?


\newcommand{\etalchar}[1]{$^{#1}$}

\end{document}